\documentclass[12pt]{amsart}
\usepackage{amsfonts,amssymb,amsthm,amsmath, hyperref}
\newcommand{\smat}[4]{\left(\begin{smallmatrix}
                 #1 & #2\\
                 #3 & #4
\end{smallmatrix}\right)}
\newcommand{\Z}{\mathbb{Z}}
\newcommand{\SL}{{\text {\rm SL}}}

\usepackage[utf8]{inputenc}
\usepackage[utf8]{inputenc}
\newtheorem{theorem}{Theorem}
\newtheorem{lemma}{Lemma}

\begin{document}

\title[Divisibility properties of coefficients of modular functions]{Divisibility properties of coefficients of modular functions in genus zero levels}
\author{Victoria Iba}
\author{Paul Jenkins}
\address{Department of Mathematics, Brigham Young University, Provo UT 84602}
\email{jenkins@math.byu.edu}
\author{Merrill Warnick}
\address{Department of Mathematics, Brigham Young University, Provo UT 84602}
\email{merrill.warnick@gmail.com}
\thanks{This work was partially supported by a grant from the Simons Foundation (\#281876 to Paul Jenkins).  The authors thank the referee for many helpful suggestions.}

\begin{abstract} We prove divisibility results for the Fourier coefficients of canonical basis elements for the spaces of weakly holomorphic modular forms of weight $0$ and levels $6, 10, 12, 18$ with poles only at the cusp at infinity.  In addition, we show that these Fourier coefficients satisfy Zagier duality in all weights, and give a general formula for the generating functions of such canonical bases for all genus zero levels.
\end{abstract}

\maketitle

\section{Introduction}

In 1949, Lehner~\cite{Lehner1, Lehner2} showed that many of the Fourier coefficients $c(n)$ of the classical $j$-function $j(z) = q^{-1} + 744 + \sum c(n) q^n$, where $q = e^{2\pi i z}$ as usual, are highly divisible by small primes dividing $n$.  For instance, he showed that $c(2^a3^b5^c7^dn) \equiv 0 \pmod{2^{3a+8}3^{2b+3}5^{c+1}7^d}$ for any positive integers $a, b, c, d$.  These divisibility results were generalized to stronger congruences modulo larger powers of these primes by Aas~\cite{Aas} and Kolberg~\cite{Kolberg1, Kolberg2}, and such congruences were shown to hold for all elements of the canonical basis $\{q^{-m} + O(q)\}$ for the space of weakly holomorphic modular functions for $\SL_2(\Z)$ by Griffin~\cite{Griffin}.

Let $N$ be a positive integer such that the group $\Gamma_0(N)$ has genus zero; thus, $N \in \{1, 2, 3, 4, 5, 6, 7, 8, 9, 10, 12, 13, 16, 18, 25\}$.  The space of modular functions of genus zero level $N$ with poles only at the cusp at $\infty$ has a similar canonical basis, and analogous divisibility results have been shown for the Fourier coefficients of the basis elements in the genus zero prime power levels $N \in \{2, 3, 4, 5, 7, 8, 9, 13, 16, 25\}$ in~\cite{Nick,Haddock, JT2,JT}.  In this paper, we prove divisibility results for the Fourier coefficients of canonical basis elements in the remaining genus zero levels $N \in \{6, 10, 12, 18\}$.  Additionally, we show that the coefficients of these basis elements satisfy Zagier duality, and we explicitly construct their generating functions.

Let $M_k^!(N)$ be the space of weakly holomorphic modular forms for $\Gamma_0(N)$ of even integer weight $k$, and let $M_k^\sharp(N)$ be the subspace of such forms with poles allowed only at the cusp at $\infty$.  If $n_0$ is the maximal order of vanishing at the cusp at $\infty$ for a modular form in $M_k^\sharp(N)$, then the space $M_k^\sharp(N)$ has a countably infinite canonical basis $\{f_{k, m}^{(N)}(z)\}_{m \geq -n_0}$ indexed by the order of vanishing $-m$ of $f_{k, m}^{(N)}(z)$ at the cusp at $\infty$.  Specifically, define $f_{k, m}^{(N)}(z)$ to be the unique form in $M_k^\sharp(N)$ with Fourier expansion given by \[f_{k, m}^{(N)}(z) = q^{-m} + \sum_{n > n_0} a_k^{(N)}(m, n) q^n,\] so that the gap between $q^{-m}$ and the next nonzero Fourier coefficient is as large as possible.  These bases were explicitly constructed for prime and prime power levels $N$ of genus zero in~\cite{Duke, Sharon, Haddock, JT2, JT}.

Let $S_k^\sharp(N)$ be the subspace of $M_k^\sharp(N)$ consisting of forms which vanish at all cusps except possibly at $\infty$, and let $n_1$ be the maximal order of vanishing at the cusp at $\infty$ for a modular form in this space.  As we did for $M_k^\sharp(N)$, we may construct a canonical basis $\{g_{k, m}^{(N)}(z)\}_{m \geq -n_1}$ for $S_k^\sharp(N)$, where $g_{k, m}^{(N)}(z)$ is the unique form in $S_k^\sharp(N)$ with Fourier expansion \[g_{k, m}^{(N)}(z) = q^{-m} + \sum_{n > n_1} b_k^{(N)}(m, n) q^n;\] again, the gap between $q^{-m}$ and the next nonzero Fourier coefficient is as large as possible.

The first theorem of this paper is that just as in prime power genus zero levels, the Fourier coefficients of these canonical bases satisfy Zagier duality.
\begin{theorem}\label{thm:duality}
For $N = 6, 10, 12, 18$, all even integers $k$, and all integers $m, n$, we have
\[a_k^{(N)}(m, n) = -b_{2-k}^{(N)}(n, m).\]
\end{theorem}
Additionally, we have the following formula for generating functions for genus zero level basis elements, generalizing the results in~\cite{Duke, Sharon, Haddock, JT2} to all genus zero levels.

\begin{theorem}\label{thm:genfn}
For every genus zero level $N$ we have
\[\sum_{m = -n_0}^\infty f_{k, m}^{(N)}(\tau)q^m = \frac{f_{k, -n_0}^{(N)}(\tau) g_{2-k, n_0+1}^{(N)}(z)}{f_{0, 1}^{(N)}(z) - f_{0, 1}^{(N)}(\tau)}.\]
\end{theorem}

The principal results of this paper give divisibility results for coefficients for weight zero basis elements for levels $6$, $10$, $12$ and $18$ by primes dividing the level.  The following congruences hold in level $6$.
\begin{theorem}\label{thm:level6}
Let $r,s$ be integers relatively prime to $2$. Then
\begin{align}
\label{cong1} &a_0^{(6)}(2^ar,2^bs) \equiv 0 \textnormal{ }(\textnormal{\emph{mod} }2^{a-b+2}) & \textnormal{\emph{ if }} a>b,\\
\label{cong2} &a_0^{(6)}(2^ar,2^bs) \equiv 0 \textnormal{ }(\textnormal{\emph{mod} }2^{2}) & \textnormal{\emph{ if }} b>a.
\end{align}
Now let $r,s$ be relatively prime to $3$. Then
\begin{align}
&a_0^{(6)}(3^ar,3^bs) \equiv 0 \textnormal{ }(\textnormal{\emph{mod} }3^{a-b+1}) & \textnormal{\emph{ if }} a>b,\\
&a_0^{(6)}(3^ar,3^bs) \equiv 0 \textnormal{ }(\textnormal{\emph{mod} }3^{1}) & \textnormal{\emph{ if }} b>a.
\end{align}
\end{theorem}

Using the $U_2$ and $U_3$ operators, these congruences in level $6$ lead to very similar congruences for levels $12$ and $18$.
\begin{theorem}\label{thm:level1218}
For $N=12,18$, let $r,s$ be integers relatively prime to $2$. If $N=18$, assume that $3|r$. Then
\begin{align}
&a_0^{(N)}(2^ar,2^bs) \equiv 0 \textnormal{ }(\textnormal{\emph{mod} }2^{a-b+2}) & \textnormal{\emph{ if }} a>b,\\
&a_0^{(N)}(2^ar,2^bs) \equiv 0 \textnormal{ }(\textnormal{\emph{mod} }2^{2}) & \textnormal{\emph{ if }} b>a.
\end{align}
Now let $r,s$ be relatively prime to $3$. If $N=12$, assume that $2|r$. Then
\begin{align}
&a_0^{(N)}(3^ar,3^bs) \equiv 0 \textnormal{ }(\textnormal{\emph{mod} }3^{a-b+1}) & \textnormal{\emph{ if }} a>b,\\
&a_0^{(N)}(3^ar,3^bs) \equiv 0 \textnormal{ }(\textnormal{\emph{mod} }3^{1}) & \textnormal{\emph{ if }} b>a.
\end{align}
If $N = 18$ and $3 \nmid r$, or if $N = 12$ and $2 \nmid r$, we have the weaker divisibility results
\begin{align}
&a_0^{(18)}(2^ar,2^bs) \equiv 0 \pmod{2^{a-b}} & \textnormal{\emph{ if }} a > b, \\
\label{weakcong}&a_0^{(12)}(3^ar,3^bs) \equiv 0 \pmod{3^{a-b}} & \textnormal{\emph{ if }} a > b.
\end{align}

\end{theorem}

Similar congruences hold for level $10$, the last genus zero level.
\begin{theorem}\label{thm:level10}
If $r,s$ are relatively prime to $2$, then
\begin{align}
\label{cong9} &a_0^{(10)}(2^ar,2^bs) \equiv 0 \textnormal{ }(\textnormal{\emph{mod} }2^{a-b+1}) & \textnormal{\emph{ if }} a>b, \\
\label{cong10} &a_0^{(10)}(2^ar,2^bs) \equiv 0 \textnormal{ }(\textnormal{\emph{mod} }2^{1}) & \textnormal{\emph{ if }} b>a.
\end{align}
If $r,s$ are relatively prime to $5$, then
\begin{align}
\label{cong11} &a_0^{(10)}(5^ar,5^bs) \equiv 0 \textnormal{ }(\textnormal{\emph{mod} }5^{a-b}) & \textnormal{\emph{ if }} a>b.
\end{align}
\end{theorem}
One obvious difference in level $10$ is that there is no $b>a$ congruence given modulo powers of $5$.  This behavior is similar to the congruences for Fourier coefficients in level $13$, as described in~\cite{JT}.

We note that because of the duality given by Theorem~\ref{thm:duality}, all divisibility results in weight $0$ immediately imply similar results for the coefficients of the weight $2$ basis elements.  In higher weights, such divisibility results are not expected due to the presence of cusp forms.

In Section~\ref{sec:basis} of this paper, we construct the canonical bases in levels $6, 10, 12$, and $18$.  In Section~\ref{sec:duality}, we prove Theorems~\ref{thm:duality} and~\ref{thm:genfn}.  Sections~\ref{sec:level6}, \ref{sec:level1218}, and \ref{sec:level10} contain the proofs of Theorems~\ref{thm:level6}, \ref{thm:level1218}, and~\ref{thm:level10}.

\section{Canonical Basis Construction}\label{sec:basis}

We now explicitly construct the canonical bases for $M_k^\sharp(N)$ and $S_k^\sharp(N)$ for $N = 6, 10, 12, 18$.  We will give explicit explanations for the construction in level $6$ and more general outlines for the other levels, as all of the constructions are similar.

The congruence subgroup $\Gamma_0(6)$ has four cusps: $\infty$, $0$, $\frac{1}{2}$ and $\frac{1}{3}$. For this level, we choose the Hauptmodul
\[ \psi^{(6)}(z) = f_{0,1}^{(6)}(z) = \frac{\eta(2z)^8\eta(3z)^4}{\eta(z)^4\eta(6z)^8}-4=q^{-1}+6q+4q^2-3q^3+O(q^4) \in M_0^\sharp(6),\]
which has a pole of order $1$ at $\infty$.
We also use the form
\[ f_{2,-2}^{(6)}(z) = \frac{\eta(z)^2\eta(6z)^{12}}{\eta(2z)^4\eta(3z)^6}=q^2-2q^3+3q^4-q^6+7q^8+O(q^9) \in M_2(6),\]
which has order of vanishing $2$ at $\infty$ and has no other zeros or poles.

For nonnegative even integers $k$, the dimension of the space $M_k(6)$ of holomorphic modular forms is $k+1$, so increasing the weight by $2$ increases the dimension by $2$.  Note that multiplying or dividing a form in $M_k^\sharp(6)$ by $f_{2, -2}^{(6)}(z)$ increases or decreases both the form's weight and order of vanishing at $\infty$ by $2$, and can only introduce poles at the cusp at $\infty$.  Thus, in weight $k$ we have the initial basis element $f_{k, -k}^{(6)}(z) = \left(f_{2, -2}^{(6)}(z)\right)^\frac{k}{2}$, which has maximal order of vanishing at $\infty$ for that space. We may recursively construct additional basis elements by multiplying by the Hauptmodul $\psi^{(6)}(z)$ and subtracting earlier basis elements to obtain the largest possible gap between the leading power of $q$ and the next nonzero coefficient.  Therefore, for all $m \geq -k$, the basis element $f_{k, m}^{(6)}(z)$ has the form
\[ f_{k,m}^{(6)}(z)= q^{-m} + \sum_{n = k+1}^\infty a_k^{(6)}(m, n) q^n = \left(f_{2,-2}^{(6)}(z)\right)^\frac{k}{2}P\left(\psi^{(6)}(z)\right), \]
where $P(x)$ is a polynomial in $x$ of degree $m+k$.  We note that because $\psi^{(6)}(z)$ and $f_{2, -2}^{(6)}(z)$ have integral Fourier coefficients, the Fourier coefficients $a_k^{(6)}(m, n)$ of these basis elements and the coefficients of the polynomial $P(x)$ are also all integers.  Additionally, we have $n_0 = k$ for the largest possible order of vanishing at $\infty$ of a form in $M_k^\sharp(6)$.

The congruence subgroup $\Gamma_0(12)$ has six cusps: $\infty$, $0$, $\frac{1}{2}$, $\frac{1}{3}$, $\frac{1}{4}$ and $\frac{1}{6}$. For this level, we choose the Hauptmodul
\[ \psi^{(12)}(z) = f_{0,1}^{(12)}(z)=\frac{\eta(4z)^4\eta(6z)^2}{\eta(2z)^2\eta(12z)^4}=q^{-1}+2q+q^3+O(q^7) \in M_0^\sharp(12), \]
which may be used to generate the basis for $M_0^\sharp(12)$.

Increasing a nonnegative even weight $k$ by $2$ increases the dimension of $M_k(12)$ from $2k$ to $2k+4$.  The form
\[
\begin{split}
f_{2,-4}^{(12)}(z)= &\frac{1}{27}\frac{\eta(z)^{10}\eta(4z)\eta(6z)^9}{\eta(2z)^7\eta(3z)^6\eta(12z)^3}+\frac{11}{72}\frac{\eta(z)^7
\eta(4z)^4\eta(6z)^9}{\eta(2z)^7\eta(3z)^5\eta(12z)^4}-\frac{1}{12}\frac{\eta(z)^4\eta(4z)^7\eta(6z)^9}{\eta(2z)^7
\eta(3z)^4\eta(12z)^5}\\
&+\frac{1}{54}\frac{\eta(z)\eta(4z)^{10}\eta(6z)^9}{\eta(2z)^7\eta(3z)^3\eta(12z)}-\frac{1}{8}\frac{\eta(z)^9
\eta(4z)^3\eta(6z)^2}{\eta(2z)^6\eta(3z)^3\eta(12z)} \\
= &q^4 + O(q^5) \in M_2(12)
\end{split}
\]
has a zero of order $4$ at $\infty$.  We compute that there are no other zeros or poles at any other cusp, and note that as eta-quotients have no poles on the upper half plane, the valence formula implies that there are no other zeros.  Thus, the first basis element of $M_k^\sharp(12)$ must be $\left(f_{2, -4}^{(12)}(z) \right)^\frac{k}{2}$.
Multiplying by $\psi^{(12)}(z)$ and row reducing as before, we construct the basis elements $f_{k, m}^{(12)}(z)$ for all $m \geq -2k$, and note that $n_0 = 2k$ in level $12$.  Again, all Fourier coefficients $a_k^{(12)}(m, n)$ are integers, because $\psi^{(12)}(z)$ and $f_{2, -4}^{(12)}(z)$ have integral Fourier coefficients by the $p$-adic Sturm bound.

The congruence subgroup $\Gamma_0(18)$ has six cusps: $\infty$, $0$, $\frac{1}{2}$, $\frac{1}{3}$, $\frac{1}{6}$ and $\frac{1}{9}$. For this level, we choose the Hauptmodul
\[ \psi^{(18)}(z) = f_{0,1}^{(18)}(z)=\frac{\eta(6z)\eta(9z)^3}{\eta(3z)\eta(18z)^3}\]
and raise or lower the weight by multiplying by the form
\[
\begin{split}
f_{2,-6}^{(18)}(z) = &\frac{25}{216}\frac{\eta(z)^8\eta(6z)^2\eta(9z)^4}{\eta(2z)^4\eta(3z)^4\eta(18z)^2}  -\frac{11}{144}\frac{\eta(z)^3\eta(6z)^8\eta(9z)^7}{\eta(2z)^3\eta(3z)^6\eta(18z)^5}-\frac{121}{972}\frac{\eta(z)^6\eta(6z)^7
\eta(9z)}{\eta(2z)^3\eta(3z)^5\eta(18z)^2}\\
 &- \frac{41}{144}\frac{\eta(z)^6\eta(6z)^2\eta(9z)^6}{\eta(2z)^3\eta(3z)^4\eta(18z)^3}+\frac{67}{144}\frac{\eta(z)^4\eta(6z)^7
 \eta(9z)^3}{\eta(2z)^2\eta(3z)^5\eta(18z)^3}\\
 & +  \frac{1}{972}\frac{\eta(2z)^9\eta(3z)^8\eta(12z)}{\eta(z)^6\eta(6z)^6\eta(9z)^2}-\frac{125}{1296}\frac{\eta(z)\eta(2z)^4
 \eta(9z)^2}{\eta(3z)\eta(6z)\eta(18z)} \\
= &q^6 + O(q^7) \in M_2(18).
\end{split}
\]
Each basis element in $M_k^\sharp(18)$ will be $\left(f_{2, -6}^{(18)}(z)\right)^\frac{k}{2}$ multiplied by a polynomial in $\psi^{(18)}(z)$, and we have $n_0 = 3k$.  Again, all Fourier coefficients $a_k^{(18)}(m, n)$ are integers.

The congruence subgroup $\Gamma_0(10)$ has four cusps: $\infty$, $0$, $\frac{1}{2}$ and $\frac{1}{5}$. For this level, we use the Hauptmodul
\[\psi^{(10)}(z) = f_{0,1}^{(10)}(z)=\frac{\eta(2z)\eta(5z)^5}{\eta(z)\eta(10z)^5}-1.\]

Increasing the weight $k$ by $4$ increases the dimension of $M_k(10)$ by $6$, with the dimension of $M_0(10)$ equal to $1$ and the dimension of $M_2(10)$ equal to $3$.  Writing $k = 4\ell + k'$ with $\ell \in \mathbb{Z}$ and $k' \in \{0, 2\}$, we use the forms
\[ f^{(10)}_{2,-2} = \frac{\eta(10z)^5}{\eta(2z)} = q^3 + O(q^4) \in M_2(10), \]
\[\begin{split}
f^{(10)}_{4,-6} = &\frac{27}{10000}\frac{\eta(2z)^{14}\eta(5z)^8}{\eta(z)^{8}\eta(10z)^{6}}-\frac{13}{4000}\frac{\eta(2z)^{11}\eta(5z)^{11}}
{\eta(z)^7\eta(10z)^7} -\frac{14}{125}\frac{\eta(z)^3\eta(2z)\eta(5z)^9}{\eta(10z)^5}\\
&+\frac{167}{625}\frac{\eta(z)^2\eta(2z)^4\eta(5z)^6}{\eta(10z)^4}+\frac{1}{625}\frac{\eta(z)^{12}\eta(5z)^4}
{\eta(2z)^6\eta(10z)^2}-\frac{7}{32}\frac{\eta(z)\eta(2z)^7\eta(5z)^3}{\eta(10z)^3}+\frac{1}{16}\frac{\eta(2z)^{10}}
{\eta(10z)^2} \\
= &q^6 + O(q^7) \in M_4(10)
\end{split}\]
to see that the first basis element in weight $k = 4\ell + k'$ will be
\[\left(f_{4, -6}^{(10)}(z)\right)^\ell \left(f_{2, -2}^{(10)}(z)\right)^\frac{k'}{2}.\]
Further basis elements are obtained, as before, by multiplying by a polynomial in $\psi^{(10)}(z)$ to get the correct Fourier expansion, with $n_0 = 6\ell + k'$ and $a_k^{(10)}(m, n) \in \mathbb{Z}$.

To construct the basis elements $g_{k, m}^{(N)}(z)$ for the spaces $S_k^\sharp(N)$, we follow the same general process of finding the basis element with the largest leading power of $q$, then multiplying by the Hauptmodul and row reducing with previous basis elements.  The first basis element $g_{k, -n_1}^{(N)}(z)$ can be constructed as
\[ g_{k, -n_1}^{(N)}(z) = f_{k, -n_0}^{(N)}(z) \prod_{\frac{a}{b}} \left(\psi^{(N)}(z) - c_{\frac{a}{b}}\right),\]
where $f_{k, -n_0}^{(N)}(z)$ is the first basis element for $M_k^\sharp(N)$, the product is over all non-$\infty$ cusps of $\Gamma_0(N)$, and $c_{\frac{a}{b}}$ is the value of $\psi^{(N)}(z)$ at the cusp $\frac{a}{b}$.  It is clear that this will be a form of the correct weight which vanishes at all non-$\infty$ cusps, and it must have maximal order of vanishing at $\infty$ in $S_k^\sharp(N)$ because $f_{k, -n_0}^{(N)}(z)$ does in $M_k^\sharp(N)$, and the existence of a Hauptmodul for $\Gamma_0(N)$ means that $f_{k, -n_0}^{(N)}(z)$ cannot vanish at any cusp other than $\infty$.  Therefore, we have $n_1 = n_0 - c(N) + 1$, where $c(N)$ is the number of cusps of $\Gamma_0(N)$.

In level $6$, the Hauptmodul has the values $0$, $1$, $9$ at the three non-infinity cusps. This may be computed using Ligozat's theorem (see~\cite{Lig}) to find that the form
\[g_{2,1}^{(6)}(z) = \frac{\eta(2z)^6\eta(3z)^8}{\eta(6z)^{10}} \in M_2^\sharp(6)\]
vanishes at all non-$\infty$ cusps, and then writing the quotient $g_{2, 1}^{(6)}(z)/f_{2, -2}^{(6)}(z)$ as a polynomial in $\psi^{(6)}(z)$. The roots of the polynomial are the values at the cusps.  Thus, the first basis element of $S_k^\sharp(6)$ will be
\[g_{k, -k+3}^{(6)}(z) = \left(f_{2, -2}^{(6)}(z)\right)^{\frac{k}{2}} \psi^{(6)}(z) (\psi^{(6)}(z)-1) (\psi^{(6)}(z) - 9).\]

In level 12, the polynomial \[(\psi^{(12)}(z))(\psi^{(12)}(z)-1)(\psi^{(12)}(z)+1)(\psi^{(12)}(z)-3)(\psi^{(12)}(z)+3)\] vanishes at all non-$\infty$ cusps, while levels 18 and 10 have the polynomials
\[ (\psi^{(18)}(z))^7-7(\psi^{(18)}(z))^4-8(\psi^{(18)}(z)),\]
\[(\psi^{(10)}(z)+2)(\psi^{(10)}(z)+1)(\psi^{(10)}(z)-3),\]
which vanish at all cusps away from $\infty$.

\section{Duality and Generating Functions}\label{sec:duality}
The proof of Zagier duality for the coefficients of these canonical bases (Theorem~\ref{thm:duality}) follows the argument in~\cite{JT2} for prime power levels, which is itself a generalization of an argument credited to Kaneko by Zagier~\cite{Zagier}.

Recall that Ramanujan's $\theta$-operator acts on a weakly holomorphic modular form $f$ by $\theta f = q \frac{d}{dq} f = \frac{1}{2\pi i}\frac{d}{dz} f$ and increases the weight by $2$, preserving modularity if $f$ is weight $0$.  With the usual slash operator $f|_k \gamma = (\det \gamma)^{k/2} (cz+d)^{-k}f\left(\frac{az+b}{cz+d}\right)$ for a matrix $\gamma = \smat{a}{b}{c}{d}$, the $\theta$-operator satisfies the intertwining relation \[\theta(f|_0 \gamma) = (\theta f)|_2 \gamma.\]
For a form $f \in M_0^\sharp(N)$, we know that $f |_0 \gamma$ is holomorphic for any matrix $\gamma \in \SL_2(\Z)$ and not in $\Gamma_0(N)$, so the intertwining relation and the fact that derivatives annihilate constants imply that $\theta$ actually maps forms in $M_0^\sharp(N)$ to forms in $S_2^\sharp(N)$.  Specifically, by examining Fourier expansions we find that
\begin{equation}\label{lemma:theta}
\theta\left(f_{0,m}^{(N)}(z)\right) =-mg_{2,m}^{(N)}(z).
\end{equation}
In weight $2$, the maximal order of vanishing $n_1$ for a form in $S_2^\sharp(N)$ is equal to $-1$ for $N = 6, 10, 12, 18$, so the functions $\theta\left(f_{0,m}^{(N)}\right) =-mg_{2,m}^{(N)}$ span the space $S_2^\sharp(N)$, and every form in $S_2^\sharp(N)$ must have no constant term in its Fourier expansion.

Now consider the product $f_{k,m}^{(N)}g_{2-k,m}^{(N)}$, which must be in the space $S_2^\sharp(N)$.  By examining the Fourier expansions of the forms in the product, we find that its constant term is $a_k^{(N)}(m,n)+b_{2-k}^{(N)}(n,m)$, which must be zero, and Theorem~\ref{thm:duality} follows.

We next prove Theorem~\ref{thm:genfn}; the proof is similar to that in~\cite{ElG}.  Fix an even integer weight $k$ and a genus zero level $N$.  Let $n_0, N_1$ be the maximal orders of vanishing at $\infty$ for forms in $M_k^\sharp(N)$, $S_{2-k}^\sharp(N)$ respectively.  Thus, we have \[f_{k, m}^{(N)}(z) = q^{-m} + \sum_{n > n_0} a_k(m, n)q^n\] for every integer $m \geq -n_0$, and
\[g_{2-k, m}^{(N)}(z) = q^{-m} + \sum_{n > N_1} b_{2-k}(m, n) q^n\] for every integer $m \geq -N_1$.
Let \[F_k(z, \tau) = \sum_{m = -n_0}^\infty f_{k, m}^{(N)}(\tau) q^m\] be the generating function for the $f_{k, m}^{(N)}(z)$.  In what follows, we often write $f_{k,m}(z)$ for $f_{k, m}^{(N)}(z)$ for simplicity.

Multiply the Fourier expansions for $f_{0, 1}(z)$ and $f_{k,m}(z)$ and write the result, which is a form in $M_k^\sharp(N)$, in terms of basis elements, remembering that the powers of $q$ less than or equal to $n_0$ completely determine the form.  Doing this, we obtain
\begin{align*}
f_{0, 1}(z) f_{k,m}(z) = &\left(q^{-1} + \sum_{n=1}^\infty a_0(1, n)q^n\right)\left(q^{-m} + \sum_{n>n_0}a_k(m, n)q^n\right) \\
= &q^{-m-1} + \sum_{n=1}^\infty a_0(1, n)q^{n-m} + \sum_{n > n_0} a_k(m, n)q^{n-1} + \sum_{n=1}^\infty \sum_{N > n_0} a_0(1, n) a_k(m, N) q^{N+n} \\
= &f_{k, m+1}(z) + a_k(m, n_0+1)f_{k, -n_0}(z) + a_0(1, n_0+m)f_{k, -n_0}(z) \\
  &+\sum_{n=1}^{n_0+m-1} a_0(1, n) f_{k, m-n}(z) \hfill
\end{align*}
for any $m \geq -n_0$, where the last sum appears only if $m \geq -n_0+2$ and the second to last term appears only if $m \geq -n_0+1$.  This gives us a recurrence relation for $f_{k, m+1}(z)$ in terms of forms $f_{k, i}(z)$ with $i \leq m$.

Using this recurrence, we may write
\begin{equation*}
\begin{split}
F_k(z, \tau) = &f_{k, -n_0}(\tau)q^{-n_0} + \sum_{m = -n_0+1}^\infty f_{k, m}(\tau) q^m = f_{k,-n_0}(z) q^{-n_0} + \sum_{m = -n_0}^\infty f_{k, m+1}(z) q^{m+1} \\
= &f_{k, -n_0}(\tau)q^{-n_0} + \sum_{m=-n_0}^\infty (f_{0, 1}(\tau) f_{k,m}(\tau)q^{m+1} - a_k(m, n_0+1)f_{k, -n_0}(\tau)q^{m+1} \\
& - a_0(1, n_0+m)f_{k, -n_0}(\tau)q^{m+1} - \sum_{n=1}^{n_0+m-1} a_0(1, n)f_{k, m-n}(\tau)q^{m+1}).
\end{split}
\end{equation*}

Using the duality $a_k(m, n_0+1) = -b_{2-k}(n_0+1, m)$, we note that the term \[-\sum_{m=-n_0}^\infty
a_k(m,n_0+1)f_{k, -n_0}(\tau)q^{m+1}\] may, if $N_1\geq -n_0-1$, be rewritten as $q \cdot f_{k, -n_0}(\tau) (g_{2-k, n_0+1}(z) - q^{-n_0-1})$.  We note also that this inequality is true (and is, in fact, an equality) for all genus zero levels.  Reversing the order of the double summation, we rewrite other terms similarly to obtain
\begin{align*}
F_k(z, \tau) = &f_{k, -n_0}(\tau)q^{-n_0} + q f_{0, 1}(\tau)F_k(z, \tau) \\
               &+f_{k, -n_0}(\tau) q (g_{2-k, n_0+1}(z) - q^{-n_0-1}) - f_{k, -n_0}(\tau) q^{1-n_0} (f_{0, 1}(z)-q^{-1}) \\
               &-q(f_{0, 1}(z) - q^{-1})(F_k(z, \tau) - f_{k, -n_0}(\tau) q^{-n_0}),
\end{align*}
which simplifies to
\[0 = qf_{0, 1}(\tau)F_k(z, \tau) + q f_{k, -n_0}(\tau)g_{2-k, n_0+1}(z)-qf_{0, 1}(z) F_k(z, \tau).\]
Therefore, we have the generating function
\[F_k(z, \tau) = \frac{f_{k, -n_0}(\tau) g_{2-k, n_0+1}(z)}{f_{0, 1}(z) - f_{0, 1}(\tau)}.\]

\section{Level 6 Congruences}\label{sec:level6}

We will now prove Theorem~\ref{thm:level6} by first proving the congruence~(\ref{cong2}), and then showing that~(\ref{cong2}) implies~(\ref{cong1}).

We use the usual slash operator defined at the beginning of section~\ref{sec:duality}, omitting the subscript $k$ for brevity.  If $f \in M_0^\sharp(6)$ and $\gamma \in \SL_2(\Z)$, this becomes $f|\gamma = f\left(\frac{az+b}{cz+d}\right)$.
Recall that if $p | N$, the $U_p$ operator
\[U_p f = \frac{1}{p} \sum_{i=0}^{p-1} f|\smat{1}{i}{0}{p}\]
sends forms in $M_k^!(N)$ to forms in $M_k^!(N)$, while if $p^2|N$, the $U_p$ operator sends forms in $M_k^!(N)$ to forms in $M_k^!(\frac{N}{p})$.  Its action on Fourier expansions is given by $U_p \sum a(n) q^n = \sum a(pn)q^n$.
Thus, the $U_2$ operator on forms in $M_0^\sharp(6)$ can be written as
\[2U_2f(z) = f|\smat{1}{0}{0}{2} + f|\smat{1}{1}{0}{2}.\]

We are interested in the order of vanishing of $U_2f(z)$ at the other cusps of $\Gamma_0(6)$, which are $0, \frac{1}{2}$, and $\frac{1}{3}$.  Computing, we find that
\[2U_2f|\smat{0}{-1}{1}{0} = f|\smat{0}{-1}{1}{0}|\smat{2}{0}{0}{1} + f|\smat{1}{0}{2}{1}|\smat{1}{-1}{0}{2}, \]
\[2U_2f|\smat{1}{0}{2}{1} = f|\smat{1}{0}{4}{1}|\smat{1}{0}{0}{2} + f|\smat{3}{-1}{4}{-1}|\smat{1}{1}{0}{2}, \]
\[2U_2f|\smat{1}{0}{3}{1} = f|\smat{1}{0}{6}{1}|\smat{1}{0}{0}{2} + f|\smat{2}{-1}{3}{-1}|\smat{2}{1}{0}{1}.\]
Under the action of $\Gamma_0(6)$, we have $\frac{1}{4}$ equivalent to $\frac{1}{2}$, and $\frac{3}{4}$ also equivalent to $\frac{1}{2}$.  Since replacing $z$ with $\frac{az+b}{d}$ does not change holomorphicity at a cusp, it follows that if $f \in M_0^!(6)$ is holomorphic at $0$ and $\frac{1}{2}$, then $U_2 f$ is as well.  Therefore, applying $U_2$ repeatedly to a form in $M_0^\sharp(6)$ will give a form in $M_0^!(6)$ with poles that may occur only at the cusps at $\infty$ and $\frac{1}{3}$.

We note that for a form in $M_0^!(6)$, slashing with the matrix $\smat{2}{-1}{6}{-2}$ is an Atkin-Lehner involution which takes a form in $M_0^!(6)$ to a form in $M_0^!(6)$.  Computing as before, we find that
\[f|\smat{2}{-1}{6}{-2}|\smat{0}{-1}{1}{0} = f|\smat{-1}{-2}{-2}{-6} = f|\smat{-1}{0}{-2}{-1}|\smat{1}{2}{0}{2},\]
\[f|\smat{2}{-1}{6}{-2}|\smat{1}{0}{2}{1} =  f|\smat{0}{-1}{2}{-2} = f|\smat{0}{1}{-1}{1}|\smat{-2}{1}{0}{-1},\]
\[f|\smat{2}{-1}{6}{-2}|\smat{1}{0}{3}{1} = f|\smat{-1}{-1}{0}{-2}.\]
These calculations show that if $f$ is holomorphic at $\infty$, then $f|\smat{2}{-1}{6}{-2}$ is holomorphic at $\frac{1}{3}$; similarly, if $f$ is holomorphic at $0, \frac{1}{2}$, or $\frac{1}{3}$, then $f|\smat{2}{-1}{6}{-2}$ is holomorphic at $\frac{1}{2}, 0$, or $\infty$ respectively.

We now show that if $b > a$ and $r, s$ are odd integers, then the coefficient $a_0^{(6)}(2^a r, 2^b s)$ is divisible by $4$, proving equation~(\ref{cong2}).  Suppose first that $a = 0$. Applying the $U_2$ operator to the form $f_{0, r}^{(6)}(z) = q^{-r} + \sum_{n \geq 1} a_0^{(6)}(r, n) q^n$ gives us a modular form in $M_0^!(6)$ with Fourier expansion $\sum_{n \geq 1} a_0^{(6)}(r, 2n) q^n$ which is clearly holomorphic at $\infty$; from the discussion above, it must also be holomorphic at $0$ and $\frac{1}{2}$, since $f_{0, r}^{(6)}$ is as well.  Thus, $U_2 f_{0, r}^{(6)}(z)$ has poles only at $\frac{1}{3}$.

Replacing $z$ with $\frac{2z-1}{6z-2}$ in the formula for $U_2 f(z)$ above, we compute that
\[2U_2 f_{0, r}^{(6)}\left(\frac{2z-1}{6z-2}\right) = f_{0, r}^{(6)}\left(z-\frac{1}{2}\right) + f_{0, r}^{(6)}\left(\frac{8z-3}{12z-4}\right).\] The first term gives $-q^{-r} + O(1)$ in the Fourier expansion since $r$ is odd. The second term is holomorphic at $\infty$ because $\smat{8}{-3}{12}{-4} = \smat{2}{-3}{3}{-4} \smat{4}{0}{0}{1}$ and $f_{0, r}^{(6)}|\gamma$ is holomorphic at $\infty$ for any $\gamma \in \SL_2(\Z)$ not in $\Gamma_0(6)$.  Slashing with $\smat{0}{-1}{1}{0}, \smat{1}{0}{2}{1}, \smat{1}{0}{3}{1}$ as above, we find that this form is holomorphic at 0 and $\frac{1}{2}$, and that at the cusp at $\frac{1}{3}$ it is equal to $f_{0, r}^{(6)}(\frac{z+1}{4}) + f_{0, r}^{(6)}(\frac{z+3}{4})$.  Replacing $z$ by $\frac{z+1}{4}$ and $\frac{z+3}{4}$ in the Fourier expansion of $f_{0, r}^{(6)}(z)$, we find that the negative powers of $q$ cancel because $r$ is odd, and conclude that the form is also holomorphic at $\frac{1}{3}$.  Since its Fourier expansion begins $-q^{-r}$, it must be $-f_{0, r}^{(6)}(z)$.

Let $\psi(z)$ be the Hauptmodul $\psi(z) = \frac{\eta(z)^5 \eta(3z)}{\eta(2z)\eta(6z)^5} \in M_0^\sharp(6)$; we note that $\psi(z)$ has a pole of order $1$ at $\infty$ and vanishes at the cusp at $0$.  From the construction of the basis elements $f_{0, m}^{(6)}(z)$, we may write each basis element as an integral polynomial in powers of $\psi(z)$, giving us
\[2U_2 f_{0, r}^{(6)}\left(\frac{2z-1}{6z-2}\right) = -f_{0, r}^{(6)}(z) = \sum_{i=0}^r c_i \psi^i(z),\]
with the coefficients $c_i \in \Z$.

Using the equalities $\eta(z+1) = e^{2\pi i/24} \eta(z)$ and $\eta(-1/z) = \sqrt{z/i}\eta(z)$, we compute that
\[\psi|\smat{2}{-1}{6}{2} = \frac{-2^3 \eta(2z)^5\eta(6z)}{\eta(3z)^5\eta(z)}.\]  We note that $\psi_{\frac{1}{3}}(z) = \frac{\eta(2z)^5\eta(6z)}{\eta(3z)^5\eta(z)}$ is a modular form in $M_0^!(6)$ with integer coefficients; it is holomorphic at $\infty$ and $0$, vanishes at $\frac{1}{2}$, and has a pole of order 1 at $\frac{1}{3}$.

Replacing $z$ with $\frac{2z-1}{6z-2}$ in our expression for $2U_2 f_{0, r}^{(6)}$ above, we find that
\[2U_2 f_{0, r}^{(6)}(z) = \sum_{i=0}^r c_i 2^{3i} \psi_{\frac{1}{3}}(z)^i.\]  Because the $c_i$ and all Fourier coefficients of $\psi_{\frac{1}{3}}$ are integers, we conclude that for $n \geq 1$, all of the Fourier coefficients $a_0^{(6)}(r, 2n)$ of $U_2 f_{0, r}^{(6)}(z)$ must be divisible by 4.

Now suppose that $a \geq 1$ and that the coefficients $a_0^{(6)}(2^\alpha r', 2^\beta s')$ are divisible by 4 for every nonnegative integer $\alpha < a$ and every $\beta > \alpha$ and odd $r', s'$.  We will show that $a_0^{(6)}(2^a r, 2^b s)$ is divisible by $4$ for every integer $b > a$, proving the theorem by induction.

Define $G(z)$ to be the form $U_2 f_{0, 2^a r}^{(6)}(z) - f_{0, 2^{a-1} r}^{(6)}(z)$.  We note that $G(z)$ is holomorphic at $\infty$ because the $q^{2^{a-1}r}$ terms cancel, that it is holomorphic at $0$ and $\frac{1}{2}$ and meromorphic at $\frac{1}{3}$, and that its $n$th Fourier coefficient is $a_0^{(6)}(2^a r, 2n) - a_0^{(6)}(2^{a-1}r, n)$.  Replacing $z$ by $\frac{2z-1}{6z-2}$ as before, we have the form
\[2 G\left(\frac{2z-1}{6z-2}\right) = f_{0, 2^a r}^{(6)}\left(z-\frac{1}{2}\right) + f_{0, 2^a r}^{(6)}\left(\frac{8z-3}{12z-4}\right) - 2 f_{0, 2^{a-1} r}^{(6)}\left(\frac{2z-1}{6z-2}\right),\]
which is holomorphic away from $\infty$.
The second two terms are holomorphic at $\infty$, and the first term gives a $q^{-2^a r}$, so this must be $f_{0, 2^a r}(z)$, which has integer coefficients.  Thus, we have
\[2 G\left(\frac{2z-1}{6z-2}\right) = \sum_{i=0}^{2^a r} c_i \psi(z)^i,\] with the coefficients $c_i \in \Z$.  Replacing $z$ with $\frac{2z-1}{6z-2}$ as before, we get
\begin{equation}2G(z) = \sum_{i=0}^{2^a r}c_i 2^{3i-1} \psi_{\frac{1}{3}}(z)^i. \label{sumeqn} \end{equation}
Applying the $U_2$ operator $a$ times to this equation, we get a form $U_2^a G(z) \in M_0^!(6)$ with poles only at the cusp at $\frac{1}{3}$ and with $n$th Fourier coefficient given by
\[a_0^{(6)}(2^a r, 2^{a+1} n) - a_0^{(6)}(2^{a-1}r, 2^a n).\]  By our inductive hypothesis, the second term in this expression is divisible by 4.  Since each positive power of $q^n$ in~\eqref{sumeqn} has a coefficient divisible by $4$ as well, we conclude that $a_0^{(6)}(2^a r, 2^{a+1} n)$ is divisible by $4$, proving congruence~(\ref{cong2}).

Now consider the coefficient $a_0(2^a r,2^b s)$ with $a > b$. By duality and equation~(\ref{lemma:theta}), we have \[  a_0(2^a r,2^b s) = \frac{-2^a r}{2^b s}b_2(2^a r,2^b s) = \frac{2^a r}{2^b s}a_0(2^b r,2^a s).\]
Note that $\frac{r}{s}a_0(2^b r, 2^a s)$ must be an integer because all basis elements have integer coefficients in level 6. Then by the congruence just proved, since $a>b$, we have $2^2|a_0(2^b r, 2^a s)$, and it follows that $2^{a-b+2}|a_0(2^a r, 2^b s)$.

The proof of congruences for powers of 3 in level $6$ is very similar; we will note explicit calculations that differ.
The $U_3$ operator on forms $f \in M_0^\sharp(6)$ can be written as
\[3U_3f(z) = f|\smat{1}{0}{0}{3} + f|\smat{1}{1}{0}{3} + f|\smat{1}{2}{0}{3}.\]
To check the order of vanishing of $U_3f(z)$ at the cusps $0$, $\frac{1}{2}$ and $\frac{1}{3}$ of $\Gamma_0(6)$, we compute that
\[3U_3f|\smat{0}{-1}{1}{0} = f|\smat{0}{-1}{1}{0}|\smat{3}{0}{0}{1} + f|\smat{1}{0}{3}{1}|\smat{1}{-1}{0}{3} + f|\smat{2}{-1}{3}{-1}|\smat{1}{1}{0}{3},\]
\[3U_3f|\smat{1}{0}{3}{1} = f|\smat{1}{0}{9}{1}|\smat{1}{0}{0}{3} + f|\smat{4}{-1}{9}{-2}|\smat{1}{1}{0}{3} + f|\smat{7}{3}{9}{4}|\smat{1}{-1}{0}{3},\]
\[3U_3f|\smat{1}{0}{2}{1} = f|\smat{1}{0}{6}{1}|\smat{1}{0}{0}{3} + f|\smat{1}{0}{2}{1}|\smat{3}{1}{0}{1} + f|\smat{5}{-1}{6}{-1}|\smat{1}{1}{0}{3}.\]

Since the cusps at $\frac{2}{3}$, $\frac{1}{9}$, $\frac{4}{9}$, and $\frac{7}{9}$ are equivalent to the cusp at $\frac{1}{3}$, if $f(z)$ is holomorphic at $0$ and $\frac{1}{3}$, then $U_3f(z)$ will be also. Thus, applying $U_3$ repeatedly to a form in $M_0^\sharp(6)$ gives a form with poles only at $\frac{1}{2}$ and $\infty$.

We again use an Atkin-Lehner involution $\smat{3}{1}{6}{3}$. Computing as before, we find that
\[f|\smat{3}{1}{6}{3}|\smat{0}{-1}{1}{0} = f|\smat{1}{-3}{3}{-6} = f|\smat{1}{0}{3}{1}|\smat{1}{-3}{0}{3},\]
\[f|\smat{3}{1}{6}{3}|\smat{1}{0}{2}{1} =  f|\smat{5}{1}{12}{3} = f|\smat{5}{2}{12}{5}|\smat{1}{-1}{0}{3} = f|\smat{1}{-1}{0}{3},\]
\[f|\smat{3}{1}{6}{3}|\smat{1}{0}{3}{1} = f|\smat{6}{1}{15}{3} = f|\smat{2}{3}{5}{8}|\smat{3}{-1}{0}{1}\]
for any $f \in M_0^!(6)$.

Suppose now that $b > a$ and that $r, s$ are not divisible by $3$.  In the $a=0$ case, replacing $z$ with $\smat{3}{1}{6}{3}$ in the formula for $U_3f(z)$, we compute that
\[3U_3 f_{0, r}^{(6)}\left(\frac{3z+1}{6z+3}\right) = f_{0, r}^{(6)}\left(z+\frac{1}{3}\right) + f_{0, r}^{(6)}\left(\frac{9z+4}{18z+9}\right) + f_{0, r}^{(6)}\left(z+\frac{2}{3}\right).\]
The second term is holomorphic at $\infty$ since $\smat{9}{4}{18}{9} = \smat{1}{2}{4}{9} \smat{9}{0}{0}{1}$. The first and last terms give $-q^{-r}+O(1)$ because
$r\not\equiv 0\ (\textrm{mod}\ 3) $. Slashing with $\smat{0}{-1}{1}{0}$ and $\smat{1}{0}{3}{1}$, we find that this form is holomorphic at $0$ and $\frac{1}{3}$. Slashing with $\smat{1}{0}{2}{1}$ gives $f_{0, r}^{(6)}(\frac{z+2}{9})+f_{0, r}^{(6)}(\frac{z+8}{9}) + f_{0, r}^{(6)}(\frac{z+5}{9})$, and all negative powers of $q$ cancel out since $r\not\equiv 0 \pmod{3}$. We conclude that $3U_3 f_{0, r}^{(6)}\left(\frac{3z+1}{6z+3}\right) =-f^{(6)}_{0,r}(z)$.

Now let $\psi(z)$ be the level $6$ Hauptmodul $\psi(z) = \frac{\eta(z)^5\eta(3z)}{\eta(2z)\eta(6z)^5}$. Because powers of $\psi$ form a basis for $M_0^\sharp(6)$, we can write
\[3U_3 f_{0, r}^{(6)}\left(\frac{3z+1}{6z+3}\right) = -f_{0, r}^{(10)}(z) = \sum_{i=0}^r c_i \psi(z)^i,\]
with the coefficients $c_i \in \Z$. We compute that
\[\psi|\smat{3}{1}{6}{3} = \frac{3 \eta(3z)^5\eta(z)}{\eta(6z)\eta(2z)^5}.\]
Note that $\psi_{\frac{1}{2}}(z) = \frac{\eta(3z)^5\eta(z)}{\eta(6z)\eta(2z)^5}$ is a modular form in $M_0^!(6)$ with integer coefficients that is holomorphic at $\infty$ and $0$, vanishes at $\frac{1}{3}$, and has a pole of order $1$ at $\frac{1}{2}$.

We then use near-identical arguments for the base case and inductive step to show that the stated congruences hold for powers of $3$.

\section{Congruences in levels $12$ and $18$}\label{sec:level1218}
To extend the level $6$ congruences to levels $12$ and $18$, we first relate basis elements in these levels to basis elements in level $6$, using the usual $U_p$ and $V_p$ operators with $p=2, 3$.  Recall that $V_p(f(z)) = f(pz)$, and that if $f \in M_k^!(N)$, then $V_p f \in M_k^!(pN)$.

\begin{lemma}\label{lemma:U_oper}
For any nonnegative integer $m$, we have
\[ U_2\left(f_{0,m}^{(12)}\right)=f_{0,\frac{m}{2}}^{(6)} \textnormal{ \emph{if} }2|m,  \]
\[ U_3\left(f_{0,m}^{(18)}\right)=f_{0,\frac{m}{3}}^{(6)} \textnormal{ \emph{if} }3|m. \]
If $m$ is odd, then $U_2(f_{0, m}^{(12)}) = 0$, and if $3 \nmid m$, then $U_3(f_{0, m}^{(18)}) = 0$.
\end{lemma}
To see this for $p=2$, let $f_{0,2m}\in M^\sharp_0(12)$. We compute that $V_2U_2$ preserves the space $M_0^\sharp(12)$ as in section 3 of \cite{Haddock}.
The first term of $f_{0, 2m}$ will be $q^{-2m}$, as will the first term of $V_2 f_{0,m}^{(6)}$. Then in the expression
\[ V_2\left(U_2\left(f_{0,2m}^{(12)}(z)\right)-f_{0,m}^{(6)}(z)\right),\]
the first terms will cancel and the resulting form will have only nonnegative powers of $q$.  Since it is in the space $M^\sharp_0(12)$, it must be identically zero.
If $2$ does not divide $m$, then $U_2 f_{0, m}^{(12)}$ will be a form in $M_0^\sharp(6)$ which vanishes at $\infty$, and must be zero; thus, when $m$ is odd, the form $f_{0, m}^{(12)}(z)$ will have only odd powers of $q$ in its Fourier expansion.  The proof of the case for level $18$ is similar.

Using Lemma~\ref{lemma:U_oper}, we can prove the congruences for levels $12$ and $18$.
Consider $f_{0,m}^{(12)}\in M_0^\sharp(12)$,
and assume that $m = 2^ar$ for some $a\geq 0$ and odd $r$. Then by Lemma~\ref{lemma:U_oper}, we have $U_2(f_{0,m}^{(12)}) = f_{0,2^{a-1}r}^{(6)}$ if $a \geq 1$ and $U_2(f_{0,m}^{(12)})=0$ if $a = 0$. If $a \geq 1$, then for some $b \geq 0$ and some odd $s$, we have the relation
\[ a_0^{(12)}(2^ar,2^bs) = a_0^{(6)}(2^{a-1}r,2^{b-1}s).\]
Note that if $b=0$, then $a_0^{(12)}(2^ar,2^bs) = 0$.
Suppose that $a>b$. If $b=0$, then the congruence holds vacuously. If $b>0$, then $a-1 > b-1$ and the congruence in level 12 follows from the congruence in level $6$. If instead $b>a$, then if $a=0$, we have $U_2 f_{0, 2^a r}^{(12)} = 0$ and the congruence vacuously holds. Otherwise, $b-1>a-1$ and the congruence again follows from the congruence in level $6$.

For the congruences modulo powers of $3$, write $m=2^a 3^c r$, with $a \geq 1$ and $r'$ not divisible by $2$ or $3$.  Using the $U_2$ operator again, we obtain the relation
\[ a_0^{(12)}(2^a3^cr,2^b3^ds) = a_0^{(6)}(2^{a-1}3^cr,2^{b-1}3^ds')\]
for appropriate $b, d, s'$ with $s'$ not divisible by $2$ and $3$. Again, the congruences follow directly from the congruences in level $6$.  The weaker congruence~(\ref{weakcong}) modulo powers of $3$ when $2$ does not divide the order of the pole can be proved in a manner similar to the way in which~(\ref{cong2}) implies~(\ref{cong1}), using duality, equation~(\ref{lemma:theta}), and the fact that all coefficients $a_0^{(12)}(m, n)$ are integers.

The proof for the level 18 congruences is similar.

\section{Congruences in Level 10}\label{sec:level10}
The proof of congruences~(\ref{cong9}) and~(\ref{cong10}) for powers of $2$ in level $10$ is very similar to the proof for level $6$ in section~\ref{sec:level6}. We will make note of explicit calculations that differ.

Let $f \in M_0^\sharp(10)$. Again, we want to find the behavior at the cusps under the $U_2$ operator.
Level $10$ has the three non-infinity cusps $0$, $\frac{1}{2}$, and $\frac{1}{5}$. The calculations showing order of vanishing for the cusps at $0$ and $\frac{1}{2}$ are similar to the level $6$ calculations. For the cusp at $\frac{1}{5}$, we find that
\[2U_2f|\smat{1}{0}{5}{1} = f|\smat{1}{0}{10}{1}|\smat{1}{0}{0}{2} + f|\smat{3}{-2}{5}{-3}|\smat{2}{1}{0}{1}. \]
It follows that applying $U_2$ repeatedly to a form in $M_0^\sharp(10)$ gives forms with poles only at $\frac{1}{5}$ and $\infty$.

For the space $M_0^!(10)$, the matrix $\smat{4}{-1}{10}{-2}$ is an Atkin-Lehner involution. Computing as before,
\[f|\smat{4}{-1}{10}{-2}|\smat{0}{-1}{1}{0} = f|\smat{-1}{-4}{-2}{-10} = f|\smat{1}{0}{2}{1}|\smat{-1}{-4}{0}{-2},\]
\[f|\smat{4}{-1}{10}{-2}|\smat{1}{0}{2}{1} =  f|\smat{2}{-1}{6}{-2} = f|\smat{1}{0}{3}{1}|\smat{2}{-1}{0}{-1},\]
\[f|\smat{4}{-1}{10}{-2}|\smat{1}{0}{5}{1} = f|\smat{-1}{-1}{0}{-2}.\]

Now if we replace $z$ with $\smat{4}{-1}{10}{-2}$ in the formula for $U_2f(z)$, we compute that
\[2U_2 f_{0, r}^{(10)}\left(\frac{4z-1}{10z-2}\right) = f_{0, r}^{(10)}\left(\frac{4z-1}{20z-4}\right) + f_{0, r}^{(10)}\left(z-\frac{1}{2}\right).\]
Note that the first term here is holomorphic at $\infty$ since $\smat{4}{-1}{20}{-4} = \smat{1}{-1}{5}{-4} \smat{4}{0}{0}{1}$ and that the second term will give $-q^{-r}+O(1)$.  After checking that this is holomorphic at the other cusps, we conclude that $2U_2 f_{0, r}^{(10)}\left(\frac{4z-1}{10z-2}\right) =-f^{(10)}_{0,r}(z)$.

Now let $\psi(z)$ be the level $10$ Hauptmodul $\psi(z) = \frac{\eta(2z)\eta(5z)^5}{\eta(z)\eta(10z)^5}$. By the construction of basis elements, we can write
\[2U_2 f_{0, r}^{(10)}\left(\frac{4z-1}{10z-2}\right) = -f_{0, r}^{(10)}(z) = \sum_{i=0}^r c_i \psi^i(z),\]
with the coefficients $c_i \in \Z$. We compute that
\[\psi|\smat{4}{-1}{10}{2} = \frac{-2^2 \eta(z)\eta(10z)^5}{\eta(2z)\eta(5z)^5}.\]
Note that $\psi_{\frac{1}{5}}(z) = \frac{\eta(z)\eta(10z)^5}{\eta(2z)\eta(5z)^5} $ is a modular form in $M_0^!(10)$ with integer coefficients that is holomorphic at $\infty$, vanishes at $\frac{1}{2}$ and $0$, and has a pole of order $1$ at $\frac{1}{5}$.

After these explicit calculations, we can use near-identical arguments for the base case and inductive step to show that congruence~(\ref{cong10}) holds.  Congruence~(\ref{cong9}) follows as before using duality and equation~(\ref{lemma:theta}); the congruence~(\ref{cong11}) modulo powers of $5$ is proved in the same way.

\bibliographystyle{amsplain}
%\bibliography{references}

\end{document}